%% file: lr_lyapunov_twicer.tex
\documentclass{siamltex}[10pt]
\usepackage{color,soul}
\usepackage{amsfonts,amsmath,amssymb,color,verbatim}
\usepackage{stmaryrd}
\usepackage{hyperref}
\usepackage[ruled,vlined]{algorithm2e}
\usepackage{tabularx,ragged2e,booktabs,caption}
\usepackage{accents}
\usepackage{pgfplots}
\usepackage{ifoddpage}
\usepackage{marginnote}
\usepackage{float} %Floats
\usepackage{placeins} %Floats!!
\usepackage[colorinlistoftodos,textwidth=4cm,shadow]{todonotes}
\usepackage{wasysym}
\usepackage{array}
\usepackage{booktabs}
\usepackage{pgf}
\usepackage{multirow}

\def\eps{\varepsilon}

\def\vec{\mathop{\mathrm{vec}}\nolimits}
\def\trace{\mathop{\mathrm{tr}}\nolimits}

\def\Span{\mathop{\mathrm{span}}\nolimits}

\newtheorem{remark}{Remark}[section]

\title{From low-rank approximation to an efficient rational Krylov subspace method for the Lyapunov equation}
\author{D. A. Kolesnikov \footnotemark[3] \and I. V.~Oseledets\footnotemark[3]\ \footnotemark[5] }

\begin{document}
\maketitle

\renewcommand{\thefootnote}{\fnsymbol{footnote}}

\footnotetext[3]{Skolkovo Institute of Science and Technology,
Novaya St.~100, Skolkovo, Odintsovsky district, 143025 
Moscow Region, Russia (denis.kolesnikov@skoltech.ru, i.oseledets@skoltech.ru)}
\footnotetext[5]{Institute of Numerical Mathematics,
Gubkina St. 8, 119333 Moscow, Russia}

\renewcommand{\thefootnote}{\arabic{footnote}}
\newcounter{Ivan}
\newcommand{\Ivan}[1]{
	\refstepcounter{Ivan}{
		\todo[inline,color={green!11!black!22},size=\small]{\textbf{[IVAN\theIvan]:}~#1}
	}
}

\newcounter{Denis}
\newcommand{\Denis}[1]{
	\refstepcounter{Denis}{
		\todo[inline,color={green!11!green!22},size=\small]{\textbf{[Denis\theDenis]:}~#1}
	}
}

\newcommand{\hlc}[2][yellow]{ {\sethlcolor{#1} \hl{#2}} }

\newcolumntype{C}[1]{>{\Centering}m{#1}}
\renewcommand\tabularxcolumn[1]{C{#1}}

\begin{abstract} 
We propose a new method for the approximate solution of the Lyapunov
equation with rank-$1$ right-hand side, which is based on extended rational
Krylov subspace approximation with adaptively computed shifts. The
shift selection is obtained from the connection between the Lyapunov
equation, solution of systems of linear ODEs and alternating least
squares method for low-rank approximation. The numerical experiments
confirm the effectiveness of our approach.
\end{abstract}

\begin{keywords}
Lyapunov equation, rational Krylov subspace, low-rank approximation, model order reduction
\end{keywords}

\begin{AMS}
65F10, 65F30, 15A06, 93B40.
\end{AMS}

\pagestyle{myheadings} \thispagestyle{plain}

\section{Introduction}

Let $A$ be an $n \times n$ stable matrix and $y_0$ is a vector of length $N$.
We consider the continuous-time Lyapunov equation with rank-$1$
right-hand side: 
\begin{equation}\label{alr:lyap}
   AX + XA^{\top} = -y_0 y^{\top}_0,
   \end{equation}
For large $n$ it is impossible to store $X$, thus a low-rank
approximation of the solution is sought:
\begin{equation}\label{alr:rank1.0}
   X \approx U Z U^{\top},  \quad U \in \mathbb{R}^{n \times r}, \quad
   Z \in \mathbb{R}^{r \times r}.
\end{equation}
Lyapunov equation has fundamental role in many application areas such
as signal processing \textcolor{black}{\cite{hinamoto-filter-1993, sanches-med-2008}} and system
and control theory 
\textcolor{black}{\cite{benner-numlineareq-2008, li-lyapappr-1999, pomet-lyapadaptreg-1992}.} 
There are many approaches for the solution of the Lyapunov equation.
Alternating directions implicit (ADI)
methods\textcolor{black}{\cite{lebedev-zolotarev-1977, lu-lyapadi-1991,
    li-lrlyapadi-2002,  jbilou-adiprec-2010}} 
are powerful techniques that arise from the solution methods for elliptic and parabolic partial
differential equations \textcolor{black}{\cite{wachspress-iterlyap-1988,
  damm-genlyapadi-2008, benner-selfshifts-2013, hochbruck-preckrylov-1995}}.

Krylov subspace methods have been successful in solving linear
systems and eigenvalues problems. They utilizes Arnoldi-type
\textcolor{black}{\cite{jbilou-proj-2006, stykel-krylov-2012, jaimoukha-krylov-1994}} or
Lanczos-type \textcolor{black}{\cite{saad-overview-1989, bai-krylov-2002}} algorithms to construct low-rank approximation using
Krylov subspaces. Krylov subspace methods have advantage in
simplicity but the convergence can be slow for ill-conditioned $A$
\textcolor{black}{\cite{simdru-lyapadapt-2009, mikkelsen-curve-2010}}. 

Rational Krylov subspace methods (extended Krylov subspace method
\cite{druskin-eks-1998,sim-krylov-2007}, adaptive rational Krylov
\cite{simdru-lyapadapt-2009,druskin-adapt-2010,druskin-rksm-2011,druskin-trks-2014}, Smith method
\cite{gugercin-modsmith-2003, penzl-cyclicsmith-1999}) are often the method
of choice. Manifold-based approaches have been proposed in
\cite{vandereycken-riemlyap-2010, vandereycken-riemopt-2010} where the solution is been sought directly in the
low-rank format \eqref{alr:rank1.0}. The main computational cost in such
algorithms is the solution of linear systems with
matrices of the form $A + \lambda_i I$. 
A comprehensive review on the solution of linear matrix equations in general and Lyapunov equation in particular can be found in \cite{simoncini-review-2013}.

In this work we start from the Lyapunov equation and a simple method
that doubles the size of $U$ at each step using the solution of an
auxiliary Sylvester equation. \textcolor{black}{The main disadvantage of this approach is that} too many linear system solvers are  required for each step.
Using the rank-$1$ approximation to the correction equation, we obtain a simple formula
for the new vector. In our experiments we also found that it is a good idea to add a Krylov vector to the subspace. This increases the accuracy significantly at almost no additional cost.  We compare the effectiveness of the new method with
the extended rational Krylov method \textcolor{black}{\cite{druskin-eks-1998, sim-krylov-2007}} and the adaptive rational Krylov
approach \textcolor{black}{\cite{druskin-adapt-2010, druskin-rksm-2011}}  on several model examples with symmetric and non-symmetric
matrices $A$ coming from discretizations of two-dimensional
\textcolor{black}{and three-dimensional} elliptic PDEs on different grids.

\section{Minimization problem}
 How do we define what is the best low-rank approximation to the Lyapunov equation? A natural way is to formulate the initial problem 
 as a minimization problem 
 $$ 
    R(X) \rightarrow \min,
 $$
and then reduce this problem to the minimization over the manifold of low-rank matrices. A popular choice is the residual:
\begin{equation}\label{alr:residual}
     R(X) = \Vert A X + X A^{\top} + y_0 y^{\top}_0 \Vert^2_F, 
 \end{equation}
 which is easy to compute for a low-rank matrix $X$. Disadvantage of
 the functional \eqref{alr:residual} is well known: it may lead to the large condition
 numbers \textcolor{black}{\cite[p. 19]{laub-schur-1979}}. For the symmetric positive definite case another functional is often used:
 $$
    R(X) =\trace(X A X) + \trace(X y_0 y^{\top}_0).
 $$
 For a non-symmetric case we can use a different functional, which is based on the connection of the low-rank solution to the Lyapunov equation
 and low-dimensional subspace approximation to the solution of a
 system of linear ODEs. Consider an ODE with the matrix $A$:
$$
   \frac{dy}{dt} = Ay, \quad y(0) = y_0.
$$
It is natural to look for the solution in the low-dimensional subspace form 
$$
y(t) \approx \widetilde{y}(t) = U c(t), 
$$
where $U$ is an $n \times r$ orthogonal matrix and $c(t)$ is an $r \times
1$ vector. Then the columns of the  \textcolor{black}{orthogonal} matrix $U$ that minimizes 
\begin{equation*}
\begin{split}
\min_{U \in O(n), c(t)} \int^{\infty}_0\limits \Vert y(t) - \widetilde{y}(t)
\Vert^2dt  &= \\
= \textcolor{black}{\min_{U \in O(n)} \int^{\infty}_0\limits \Vert y(t) - UU^{\top} y(t)
\Vert^2dt} &+ \textcolor{black}{\min_{U \in O(n), c(t)} \int^{\infty}_0\limits \Vert
 U U^{\top} y(t) - \widetilde{y}(t) \Vert^2dt} = \\
 = \min_{U \in O(n)} \int^{\infty}_0\limits \Vert y(t) - UU^{\top}
y(t) \Vert^2dt &=  \textcolor{black}{\min_{U \in O(n)} \Big( \trace X  - \trace U^{\top} X U \Big)}
\end{split}
\end{equation*}
are the eigenvectors, corresponding to the largest eigenvalues of the matrix $X$\textcolor{black}{\cite{sorensen-sylv-2002}} that solves
the Lyapunov equation \eqref{alr:lyap}. However, the computation of
\textcolor{black}{ the optimal $\widehat{c}(t) = U^{\top}y(t)$} requires the knowledge of the true solution $y(t)$ which
is not known. Instead, we can consider the Galerkin projection:
\begin{equation}\label{alr:odeappr}
\begin{split} 
\widetilde{y} &= U c(t),\\
 \frac{dc}{dt} &= U^{\top} A U c, \quad c_0 = U^{\top} y_0, 
\end{split}
\end{equation}
The final approximation is then
\begin{equation}\label{alr:odefinappr}
\begin{split}  
    \widetilde{y} = U e^{Bt} U^{\top} y_0,
\end{split}
\end{equation}
 where $B = U^{\top} A U$.
 The functional to be minimized is
 \begin{equation}\label{alr:fun1}
 F(U) = \int^{\infty}_0\limits\Vert y - \widetilde{y} \Vert^2 dt.
\end{equation}
Note that the functional depends only on $U$. Given $U$, the
approximation to the solution of the Lyapunov equation can be recovered
from the solution of the
``small''  Lyapunov equation
\begin{equation} \label{lyap:appr} 
\begin{split}
X \approx U Z U^{\top}, \quad
     B Z + Z B^{\top} = -c_0 c^{\top}_0.
\end{split}
\end{equation}
The functional \eqref{alr:fun1} can not be efficiently computed. However, a simple expansion of the norm gives
 \begin{equation}\label{alr:FU}
    F(U) =  \int^{\infty}_0\limits\Vert y \Vert^2 dt - 2 \int^{\infty}_0\limits\langle y, \widetilde{y} \rangle dt + \int^{\infty}_0\limits \Vert \widetilde{y} \Vert^2 dt. 
\end{equation}
%The first term in \eqref{alr:FU} does not depend on $U$, so it can be omitted in the minimization. 
\textcolor{black}{
Let us represent the functional in the next form}
\begin{equation*}
    F(U) = F_1(U) - 2 F_2(U),  
\end{equation*}
where
$$
F_1(U) =\int^{\infty}_0\limits (\Vert y \Vert^2 - \Vert \widetilde{y}
\Vert^2) dt  , \quad F_2(U) = \int^{\infty}_0\limits(\langle y,
\widetilde{y} \rangle -  \Vert \widetilde{y} \Vert^2) dt.
$$
\begin{lemma} \textcolor{black}{Assume that matrices A and B are stable.}
Then the functionals $F_1(U), F_2(U)$ can be calculated as follows:
\begin{equation}\label{f12}
\begin{split}
F_1(U) &= \trace X - \trace Z,\\
F_2(U) &= \trace U^{\top} (P - U Z),
\end{split}
\end{equation}
where $P$ is the solution of the Sylvester equation and $Z$ is solution of
the Lyapunov equation:
\begin{equation}\label{auxiliary1}
\begin{split}
A P + P B^{\top} &= -y_0 c_0^{\top},\\
B Z + Z B^{\top} &= -c_0 c_0^{\top}.
\end{split}
\end{equation}
\end{lemma}
{\it Proof: \\}
It is easy to see, that
\begin{equation*}
\begin{split}
 -\int^{\infty}_0\limits\langle y, \widetilde{y} \rangle dt
&= -\trace \int^{\infty}_0\limits e^{At}y_0 c_0^{\top}
e^{B^{\top}t}U^{\top} dt =\left. -\trace \left(e^{At}P e^{B^{\top}t}
  U^{\top} \right)\right|^{\infty}_0 = \trace U^{\top} P,\\
 \int^{\infty}_0\limits \Vert \widetilde{y} \Vert^2 dt &= 
\trace \int^{\infty}_0\limits U e^{Bt}c_0 c_0^{\top}
e^{B^{\top}t}U^{\top} dt = -\trace \left. \left(e^{Bt}Z 
e^{B^{\top}t}\right)\right|^{\infty}_0 = \trace Z. 
\end{split}
\end{equation*}
In the same way $\int^{\infty}_0\limits \Vert y \Vert^2 dt =
\trace X$ and we can write
\begin{equation*}
\begin{split}
F_1(U) &= \int^{\infty}_0\limits \Big(\Vert y \Vert^2 - \Vert \widetilde{y}
\Vert^2\Big) dt = \trace X - \trace Z,\\
F_2(U) &= \int^{\infty}_0\limits\Big(\langle y,
\widetilde{y} \rangle -  \Vert \widetilde{y} \Vert^2\Big) dt = \trace
(U^{\top}P - Z) =\trace U^{\top} (P - U Z). \blacksquare
\end{split}
\end{equation*}
\textcolor{black}{Note that using integral representation of the solution of Lyapunov
equation is well-known fact \cite{saad-numsol-1989}:}
\begin{equation*}
\begin{split}
X = -\int_0^{\infty} \limits e^{At}y_0 y_0^{\top} e^{A^\top t} dt.  
\end{split}
\end{equation*}
\textcolor{black}{The Sylvester equation can be solved using a standard method
\cite[Theorem 3.1]{sorensen-sylv-2002}} since the
matrix
$B$ is $r \times r, r \ll n$. We compute the Schur decomposition of $B^{\top}$ and the
equation is reduced to $r$ linear systems with the matrices $A + \lambda_i
I, 
\quad i = 1, \ldots, r.$
\begin{lemma} \textcolor{black}{Assume that matrices A and B are stable.}
The gradient of $F(U)$ can be computed as:
\begin{equation}\label{f:grad}
\begin{split}
\grad F(U) = -2 P + 2 y_0 (
c_0^{\top} Z_I - y_0^{\top} P_U) ) &+ 2 A U (Z Z_I  - P^{\top} P_U) + \\
&+ 2 A^{\top} U ( Z_I Z - P_U^{\top} P),
\end{split}
\end{equation}
where $P, Z$ are defined by \eqref{auxiliary1} and
\begin{equation}\label{auxiliary2}
\begin{split}
A^{\top} P_U + P_U B &= -U, \\
B^{\top} Z_I + Z_I B &= -I_r.
\end{split}
\end{equation}
\end{lemma}
{\it Proof: \\}
Variation of Z can be expressed as a solution of the Lyapunov equation with
another right hand side:
$$ B\delta Z + \delta Z B^{\top} = - \delta c_0 c_0^{\top} - c_0 
\delta c_0^{\top} - \delta B Z - Z \delta B^{\top}   
$$ 
Using the well-known integral form of the solution of the Lyapunov equation we get that:
\begin{equation*}
\begin{split}
-\trace\delta Z &= -\trace \int_0^{\infty} e^{Bt} (\delta c_0 c_0^{\top} + c_0
\delta c_0^{\top} + \delta B Z + Z \delta B^{\top}) e^{B^{\top} t} dt
=  \\
 &= -\trace (\int_0^{\infty} e^{B^{\top} t} I e^{Bt} dt) (\delta c_0 c_0^{\top} + c_0
\delta c_0^{\top} + \delta B Z + Z \delta B^{\top}) = \\
 &= -\trace Z_I (\delta c_0 c_0^{\top} + c_0
\delta c_0^{\top} + \delta B Z + Z \delta B^{\top}) = \\
&= -2 \trace\delta U^{\top} (y_0 c_0^{\top} Z_I + A U Z Z_I + A^{\top} U Z_IZ). 
\end{split}
\end{equation*}
Similarly for $P$:
$$ A \delta P + \delta P B^{\top} = -y_0 \delta c_0^{\top} - P \delta
B^{\top}, $$
therefore,
\begin{equation*}
\begin{split}
\delta \trace (U^{\top} P) &=\trace (\delta U^{\top} P + U^{\top}
\delta P ) = \\
&= \trace (\delta U^{\top} P + U^{\top} \int_0^{\infty} e^{At} (y_0 \delta c_0^{\top} + P \delta
B^{\top}) e^{B^{\top}t}dt) = \\
 &=\trace \delta U^{\top} P + \trace (\int_0^{\infty} e^{B^{\top}t} U^{\top} e^{At} dt)(y_0 \delta c_0^{\top} + P \delta
B^{\top}) = \\
 &= \trace \delta U^{\top} P + \trace P_U^{\top} (y_0 \delta c_0^{\top} 
+ P \delta B^{\top}) = \\
&=\trace \delta U^{\top} (P + y_0 y_0^{\top} P_U + A^{\top} U P_U^{\top} X + A U P^{\top} P_U).
\end{split}
\end{equation*}
Finally,
\begin{equation*}
\begin{split}
\delta F(U) = \trace \delta U^{\top} (\grad F) = \trace
\delta(-2 U^{\top} P &+ Z),\\
\grad F(U) = -2 P + 2 y_0 (
c_0^{\top} Z_I - y_0^{\top} P_U) ) &+ 2 A U (Z Z_I  - P^{\top} P_U) + \\
&+ 2 A^{\top} U ( Z_I Z - P_U^{\top} P). \blacksquare
\end{split}
\end{equation*}
Now denote by $R_1(U)$ and $R_2(U)$ residuals of the Lyapunov and
Sylvester equations: 
\begin{equation}\label{residual}
\begin{split}
R_1(U) &= \Vert A (U Z U^{\top}) + (U Z U^{\top}) A^{\top} + y_0
y_0^{\top}\Vert, \\
R_2(U) &= \Vert A (UZ) + (UZ) B^{\top} + y_0 c_0^{\top} \Vert.
\end{split}
\end{equation}
\begin{lemma}\label{alr:lemres}
Assume that \textcolor{black}{matrices A and B are stable and } $y_0$ lies in the column space of  $U$.
 Then the next equality holds:
\begin{equation}\label{residual:equality}
R_1(U) = \sqrt{2} R_2(U)
= \sqrt{2} \Vert (A U - U B) Z\Vert .
\end{equation}
\end{lemma}
{\it Proof: \\}
Since $Z$ is the \textcolor{black}{unique} solution of the Lyapunov equation, we get that
\begin{equation*}
\begin{split}
R_1(U)^2 &= \Vert A U Z U^{\top} + U Z U^{\top} A^{\top} + y_0
y_0^{\top}\Vert^2 = \\ 
&=\Vert y_0 y_0^{\top} - U U^{\top} y_0 y_0^{\top} U
U^{\top} + (A U - U B)Z U^{\top} + U Z (A U - U B)^{\top} \Vert^2 = \\
&= \Vert (A U - U B)Z U^{\top} \Vert^2 + \Vert U Z (A U - U B)^{\top}
\Vert^2 = 2 \Vert  (A U - U B)Z\Vert^2.
\end{split}
\end{equation*}
We can use the same trick for the residual of the Sylvester equation:
\begin{equation*}
\begin{split}
R_2(U)^2 &= \Vert A (UZ) + (UZ) B^{\top} + y_0 c_0^{\top} \Vert^2 =\\
&= \Vert (I - U U^{\top})y_0 c_0^{\top} + (A U - U B) Z\Vert^2 = \Vert
(A U - U B) Z\Vert^2. \blacksquare
\end{split}
\end{equation*}
\textcolor{black}{Lemma \ref{alr:lemres} is well known for the special case when $U$ is the basis of a (rational) Krylov
subspace \cite[Theorem 2.1]{jaimoukha-krylov-1994}.}
Lemma \ref{alr:lemres} is valid if $y_0 \in \Span U$ and we will always make sure  that  $y_0 =
U U^{\top} y_0$. 
The next Lemma shows that if the residual of the Lyapunov equation goes to zero,
so does the values of the functional $F(U)$.
\begin{lemma}\label{alr:lemmabound}
Assume that \textcolor{black}{ matrices A and B are stable and} $y_0$ lies in the column space of  $U$.
 Then $$ F(U) \leq C R_1(U) $$
 with a constant 
\textcolor{black}{
$$C = \Vert \Big(I \otimes A + A \otimes I \Big)^{-1}\Vert_F +
\sqrt{2r} \Vert \Big(I \otimes A + B \otimes I \Big)^{-1}\Vert_F .$$ 
}
%that depends only on $A$ and $r$.
\end{lemma}

{\it Proof: \\}
The matrix $X - U Z U^{\top}$ satisfies the Lyapunov equation
\begin{equation*}
A (X - UZU^{\top}) + (X - UZU^{\top}) A^{\top} = -(A U - U B)Z U^{\top}
-  U Z (A U - U B)^{\top},
\end{equation*}
therefore
$$
(I \otimes A + A \otimes I) \vec \Big(X - UZU^{\top}\Big) = -\vec \Big( (A U - U B)Z U^{\top}
+  U Z (A U - U B)^{\top} \Big), 
$$
and 
$$
\Vert X - U Z U^{\top} \Vert_F \leq \Vert \Big(I \otimes A + A \otimes I \Big)^{-1} \Vert_F ~R_1(U) 
$$
Thus, 
$$| \trace X - \trace Z | \leq C_1 R_1(U).$$
In the same way, 
$$
   \Vert P - U Z \Vert_F \leq \Vert \Big(I \otimes A + B \otimes I \Big)^{-1}\Vert_F~R_2(U) = C_2 R_2(U),
$$
therefore
$$
|F_2(U)| = |\trace U^{\top} (P - UZ)| \leq \Vert U \Vert_F ~\Vert P - U Z \Vert_F \leq \sqrt{r} C_2 R_2(U). 
$$
Finally, 
$$
   F(U) \leq |F_1(U)| + 2|F_2(U)| \leq (C_1 + \sqrt{2r} C_2) R_1(U)
   \textcolor{black}{ = C R_1(U)}.
$$
% \Vert X - UZU^{\top} \Vert &\leq \Vert \Big(I \otimes A + A \otimes
% I\Big)^{-1}\Vert R_1(U), \\
%F_1(U) = \trace X - \trace Z &\leq C_1 \Vert X - UZU^{\top} \Vert,  
%
$\blacksquare$

Lemma~\ref{alr:lemmabound} shows that $R_1(U)$, the residual of the Lyapunov equation, is a viable
error bound for our functional $F(U)$. 

%We show below that functional $F(U)$ is measure of matrix X approximation accuracy. We need
%auxiliary lemma for this purpose:

%\begin{lemma}\label{alr:fnorm}
%Let us consider $n \times r$ matrices $V$ and $W$. Assume that matrix
%$W$ has rank-$1$ structure $W = w q^{\top}.$ Then the next
%inequality holds: 
%$$
%\sqrt{2} \Vert V \Vert_F \Vert W \Vert_F \leq \Vert V W^{\top} + W V^{\top} \Vert_F
%$$
%\end{lemma}
%{\it Proof: \\}
%\begin{equation*}%\label{alr:fbound1}
%\begin{split}
%&\Vert V W^{\top} + W V^{\top} \Vert_F^2 = \\ 
%&= \sum_{1\leq i, j \leq n} \Big( \sum_{k =1}^r V_{ik} W_{jk} + V_{jk} W_{ik} \Big)^2 = 
% \sum_{1\leq i, j \leq n} \Big(\sum_{k_1 = 1}^r (V_{ik_1}
% W_{jk_1} ) + \sum_{k_2=1}^r (V_{jk_2} W_{ik_2})\Big)^2 = \\
%&= \sum_{1\leq i, j \leq n} \Big[ \Big(\sum_{k_1 = 1}^r V_{ik_1}
% W_{jk_1}  \Big)^2 + \Big(  \sum_{k_2 = 1}^r V_{jk_2}
% W_{ik_2}\Big)^2 + 2 \Big(\sum_{k_1 = 1}^r V_{ik_1}
% W_{jk_1}  \Big) \Big( \sum_{k_2 = 1}^r V_{jk_2}
% W_{ik_2} \Big) \Big] = \\
% &=  2\sum_{1\leq i, j \leq n} \Big[ \Big(\sum_{k = 1}^r V_{ik}
% W_{jk}  \Big)^2 + \sum_{k_1 = 1}^r  \sum_{k_2 = 1}^r \Big( V_{ik_1}
% w_j q_{k_1}\Big) \Big( V_{jk_2} w_i q_{k_2}\Big)\Big] = \\
%&= 2\sum_{1\leq i, j \leq n}  \Big(\sum_{k = 1}^r V_{ik}
% W_{jk}  \Big)^2 + \sum_{i = 1}^n\sum_{k = 1}^r \Big( w_i V_{ik}
% q_k\Big)^2 \geq  2\sum_{1\leq i, j \leq n}  \Big(\sum_{k = 1}^r V_{ik}
% W_{jk}  \Big)^2 = 2 \Vert V W^{\top} \Vert_F^2
%2\sum_{1\leq i, j \leq n}  \Big(\sum_{k = 1}^r V_{ik}  w_j q_k \Big)^2 = 2 \Vert V \Vert_F \Vert W \Vert_F
%\end{split}
%\end{equation*}
%$\blacksquare$
\begin{lemma}\label{alr:funcbound}
    \textcolor{black}{
Assume that matrices A and B are stable, $y_0$ lies in the column
space of  $U$, matrices $P$ and $Z$ are solutions of the Sylvester and
the small Lyapunov equations \ref{auxiliary1} correspondingly. Then $$\Vert X - P U^{\top} - U P^{\top} +
U Z U ^{\top} \Vert_F \leq F(U).$$}
\end{lemma}
{\it Proof: \\}
\textcolor{black}{
Let us consider the matrix 
$$M = \int_0^{\infty}\limits \Big(y(t)-\widetilde y(t)\Big)
\Big(y^{\top}(t)-\widetilde y^{\top}(t)\Big) dt.$$
We use exponential representations of the Lyapunov and the Sylvester
equations solutions \eqref{auxiliary1}:
\begin{equation*}%\label{alr:fbound1}
\begin{split}
M &= \int_0^{\infty}\limits y(t) y^{\top}(t) dt
-\int_0^{\infty}\limits\widetilde y(t) y^{\top}(t) dt
-\int_0^{\infty}\limits y(t) \widetilde y^{\top}(t) dt 
+ \int_0^{\infty}\limits \widetilde y(t) \widetilde y^{\top}(t) dt =
\\
&= X - P U^{\top} - U P^{\top} +
U Z U ^{\top}.
\end{split}
\end{equation*}
Matrix $M$ is symmetric and non-negative definite and therefore:
\begin{equation}\label{alr:fbound1}
\begin{split}
\Vert X - P U^{\top} - U P^{\top} +
U Z U ^{\top} \Vert_F = \Vert M \Vert_F \leq \trace M &= F(U)
\end{split}
\end{equation}
}
$\blacksquare$

\textcolor{black}{
We can also compare $F(U)$ functional based on Galerkin projection with
functional based on optimal projection $\widehat{y} = U U^{\top} y(t)$:  
\begin{equation}\label{alr:fbound2}
\begin{split}
\int_0^{\infty}\limits \Vert y(t)-\widehat y(t) \Vert^2 dt &= 
%\trace X - \trace U^{\top} X U 
\trace (I - U U^{\top}) X (I - U U^{\top}) = \\
 &= \trace (I - U U^{\top}) M (I - U U^{\top}) \leq \trace M = F(U).
\end{split}
\end{equation}
Lemma~\ref{alr:funcbound} shows that if $F(U)$ is small, the solution $X$ can be well-approximated 
by a rank-$2r$ matrix.}  
\section{Methods for basis enrichment}
Notice that the column vectors of the gradient are a linear combination of
the column vectors of $P, AU, A^{\top}U$ and $y_0.$ We need a method to
enlarge the basis U so the first idea is to use matrix $P$ to extend the basis. 
Note, that the matrix $UZ$ can be also considered as an approximation to the solution of
Sylvester equation. So to enrich the basis we will use $P_1 = P - UZ$
instead, and the matrix $P_1$ satisfies the equation: 
\begin{equation}\label{sylvester:lr}
\begin{split}
& A(P - U Z) + (P - U Z) B^{\top} = -y_0 c_0^{\top} - A U Z - U
ZB^{\top} = \\
 &= - (I - U U^{\top}) y_0 c_0^{\top} - (A U -  U B) Z, \\
&A P_1 + P_1 B^{\top} = - (A U - U B) Z,
\end{split}
\end{equation}
where we have used that $y_0 = U U^{\top} y_0$.\\
The method is summarized in Algorithm \ref{twicer:int}.
The main disadvantage of Algorithm~\ref{twicer:int} is that the
computational cost grows dramatically at each step. If $U$ has $r$ columns, 
the next step will require $r$ solutions of $n\times n$ linear systems with
matrices of the form $A + \lambda_i I$.

\begin{algorithm}[H] \label{twicer:int}
\DontPrintSemicolon
\SetKwComment{tcp}{$\triangleright$~}{}
\KwData{Input matrix $A \in \mathbb{R}^{n \times n}$, vector $y_0 \in \mathbb{R}^{n \times 1}$, maximal rank $r_{\max}$, accuracy parameter $\eps$.}
\KwResult{Orthonormal matrix $U \in \mathbb{R}^{n \times r}$.}
\Begin{
\nl set $U = \frac{y_0}{\Vert y_0 \Vert}$ \tcp*[r]{Initialization}

\nl \For{$\rank U \leq r_{\max}$}
{ 
\nl Compute $c_0 = U^{\top}y_0 ,\quad B = U^{\top}AU$

\nl Compute $Z$ as Lyapunov equation solution: $B Z + Z B^{\top} = -
c_0 c_0^{\top}$

\nl Compute error estimate $\delta = \Vert (A U - U B) Z \Vert$

\nl \If{$\delta\leq \eps$} 
{
\nl Stop
}

\nl Compute $P_1$ as Sylvester equation solution: $A P_1 + P_1 B^{\top} = -(A U - U B) Z$

\nl Update $U = \mathrm{orth}[U, P_1]$
}
}
\caption{The doubling method}
\end{algorithm}
\textcolor{black}{
    \begin{remark}
    Algorithm \ref{twicer:int}
is similar to the IRKA method \cite{beattie-irkah2-2007,
  gugercin-irka-2005, gugercin-irkareduction-2008,
  flagg-irkaconv-2012}. The IRKA method starts from initial orthogonal matrices
$V_0, W_0$ with the given rank $r$ and replaces the matrices $V_{i-1}$ and
$V_{i-1}$ by the new ones $V_i$ and $W_i$ that are the 
union of the rational
Krylov subspaces $(A + s_i I)^{-1}b$ and $(A^{\top} + s_i I)^{-1}c$, correspondingly, in every step. 
The shifts $s_i$ are obtained as eigenvalues
of the matrix $V_{i-1}^{\top} A W_{i-1}$. The main difference of Algorithm
\ref{twicer:int} algorithm from IRKA is that IRKA is usually used with the fixed rank. 
\end{remark}
}

To reduce the number of solvers required by the algorithm, we propose
two improvements. The first is to add the last Krylov vector to the subspace.
In this case, as we will show, the residual will always have rank-$1$. 
The second improvement is to add only one vector each time. In order to do 
so, we will use a simple rank-$1$ approximation to $P_1$.
\subsection{Adding a Krylov vector and a rational Krylov vector to the
  subspace}
\textcolor{black}{
In the following section we will show, that under a special basis enrichment
strategy, the rank-$1$ approximation to $P_1$ can be replaced by adding 
one vector of the form $(A + sI)^{-1} w$ to the subspace. 
In this case, as it was already described in \cite{sim-krylov-2007, simdru-lyapadapt-2009}, adding a Krylov vector preserves the 
rank-$1$ structure of the residual of equation \eqref{sylvester:lr}.  
Next Lemma is a generalization of the well-known
fact of rank-$1$ residual in the Arnoldi iteration
\cite{saad-methodsbook-2003}.
}
\begin{lemma}\label{lyap:1rkres}
Let $A$ be an $n\times n$ \textcolor{black}{stable} matrix and $y_0$ is a vector of size $n$.
Assume that an $n\times r$ orthogonal matrix $U$ and vectors $w, v$
 of size $n$ satisfy the following equations:
 $$ (I-UU^{\top}) A U =w q^{\top} , \quad v = (A + s I)^{-1}w.$$
Let us denote by $U_1$ and $U_2$ the basis of the spans of the columns of the matrices $\left[U,
  v\right]$ and $\left[U, v, w\right]$ accordingly. Then the following equality holds
  $$ \left(I-U_2U_2^{\top}\right) A U_2 = \left(I-U_2U_2^{\top}\right)Aw \widehat{q} ^{\top}.$$
\end{lemma}
{\it Proof: \\}
Due to the fact that  $(I-UU^{\top})AU = w q^{\top}$ we get that $(I-U_2U_2^{\top})AU = 0.$\\ 
On the other hand we have $$ (I-U_2U_2^{\top})Av = (I-U_2U_2^{\top})((A+sI)v -s v) = (I-U_2U_2^{\top})(w - s v) = 0.$$
Therefore $(I-U_2U_2^{\top})A U_2 = (I-U_2U_2^{\top})Awq^{\top}. $
\textcolor{black}{Note that if $U_1$ and $U_2$ are obtained by using the Gram-Schmidt process
then \\$(I-U_2U_2^{\top})A U_1 = 0$ and the matrix $(I-U_2U_2^{\top})A U_2$
has only one non-zero column and this column is the last one.}
$\blacksquare$\\
\textcolor{black}{The statements about residual rank-structure
are well known \cite{jaimoukha-krylov-1994, sim-krylov-2007}.}
Lemma \ref{lyap:1rkres} shows that if the approximation algorithm starts from $U_0 =
\frac{y_0}{\Vert y_0 \Vert}$ and adds a vector from the Krylov
subspace and a corresponding vector from the rational Krylov subspace at
each step then the residual matrix $(AU - UB)Z = (I - UU^{\top})AUZ$
is rank-$1$ at any step. 
That is the main benefit from using the
Krylov subspaces in our approach. Moreover, the residual has the form 
$\widehat{w} \widehat{q}^{\top}$, where $\widehat{w} = (I - U_1 U^{\top}_1) A w$
is the next \emph{Krylov vector}. This fact is important and will be
used later.

\subsection{Rank-$1$ approximation to the correction equation}
Suppose we have the matrix $U$ constructed by sequential addition of 
rational Krylov and Krylov vectors to the subspace and the equation for $P_1$ has the form \eqref{sylvester:lr}
Note, that due to the equality 
$$
   (A U - UB) = (I - UU^{\top}) AU,
$$
the matrix $(AU - UB) = w q^{\top}$ has only one non-zero column. 
If the new vectors are added as the rightmost vectors, 
$$
(AU - UB) Z = w z^{\top},
$$
where $z^{\top}$ is the last row of the matrix $Z$.  
Therefore, the equation for $P_1$ takes the form
\begin{equation}\label{alr:rank1corr}
AP_1 + P_1 B^{\top} = -(AU - UB) Z = w z^{\top}.
\end{equation}
The last step from the doubling method to the final algorithm is find to a
rank-$1$ approximation to the solution of \eqref{alr:rank1corr}. If $U$ is known,  
we apply one step of alternating iterations,  looking for the solution 
in the form $P_1 \approx v q^{\top}$, where $q = \frac{z}{\Vert z \Vert}$ is the
normalized last row of the matrix $Z$. \textcolor{black}{This choice of $q$ is natural
due to the right hand side of \eqref{alr:rank1corr}: $-(AU - UB) Z = \Vert z \Vert w q^{\top}.$}
The Galerkin condition for $v$ leads to equation
\begin{equation*}
    \left(A + (q^{\top} B^{\top} q) ~I\right) v= w.
\end{equation*}
\textcolor{black}{
This is the main formula for the next shift.}
Due to the simple rank-$1$ structure of the residual its norm can be
efficiently computed as 
$$R_1(U) = \sqrt{2} \Vert w \Vert ~\Vert  z \Vert. $$
The final algorithm which we call \emph{alternating low rank (ALR) method} is
presented in Algorithm \ref{al:int}. The main computational cost is the solution  of the shifted linear system.
\begin{algorithm}[H] \label{al:int}
\DontPrintSemicolon
\SetKwComment{tcp}{$\triangleright$~}{}
\KwData{Input matrix $A \in \mathbb{R}^{n \times n}$, vector $y_0 \in \mathbb{R}^{n \times 1}$, maximal rank $r_{\max}$, accuracy parameter $\eps$.}
\KwResult{Orthonormal matrix $U \in \mathbb{R}^{n \times r}$.}
\Begin{
\nl set $U = \frac{y_0}{\Vert y_0 \Vert}, w_0 = \frac{y_0}{\Vert y_0 \Vert}$ \tcp*[r]{Initialization}

\nl \For{$ \rank U \leq r_{\max}$}
{ 
\nl Compute $w_k = (I_n - UU^{\top})Aw_{k-1}$

\nl Compute $c_0 = U^{\top}y_0 ,\quad B = U^{\top}AU$

\nl Compute $Z$ as Lyapunov equation solution: $B Z + Z B^{\top} = -c_0 c_0^{\top}$

\nl Compute $z$ as the last row of the matrix $Z$.

\nl Compute error estimate $\delta = \Vert w_k\Vert ~\Vert z\Vert$

\nl \If{$\delta \leq \eps$} 
{
\nl Stop
}

\nl Compute shift $s = q^{\top} B q, q = \frac{z}{\Vert z \Vert}$

\nl Compute $v_k = (A + s I)^{-1}w_k$

\nl Update $U := \mathrm{orth}[U, v_k, w_k].$
}
}
\caption{The Adaptive low-rank method}
\end{algorithm}

\textcolor{black}{ 
\textbf{Remark.} Restriction to rank-$1$ right hand sides is helpful because there
are several ways to generalize ALR algorithm for an arbitrary rank case. 
We can either go for rank-$1$ approximation to the solution of the Sylvester equation (and get one new Rational Krylov vector each time),
or use a block version of the ALR algorithm.  We plan to compare this
approaches in our future work.
}

\section{Numerical experiments}
We have implemented the ALR method in Python using SciPy and NumPy
packages available in the Anaconda Python distribution. The
implementation is available online at 
\url{https://github.com/dkolesnikov/rkm_lyap}. The
matrices, Python code and IPython notebooks which reproduce all the figures in this work are available at \url{https://github.com/dkolesnikov/alr-paper}, where the .mat files with test matrices and vectors can be found as well.
We have compared the ALR method with two methods with publicly
available \textcolor{black}{Matlab implementations. We have ported all the methods to Python for a fair comparison, 
and their implementations are also online}.

The first method uses the Extended Krylov Subspaces approach which was
proposed in \cite{druskin-eks-1998}.
 Its main idea is to use as the basis the \emph{extended Krylov
   subspace} of the form 
$$\Span \Big(A^{-k} y_0, \ldots, y_0, \ldots,
 A^l y_0 \Big).$$ 
Note, that  the residual in this approach is also rank-$1$
 and can be cheaply computed. This approach was implemented as the 
Krylov plus Inverted Krylov algorithm (here-after KPIK) in
\cite{sim-krylov-2007} and convergence estimate also was obtained.

The second approach is the Rational Krylov Subspace Method (RKSM) which was
proposed in \cite{druskin-adapt-2010}. Its main idea is to compute
vectors step by step from the rational Krylov subspaces 
$$\Span \Big((A + s_i I)^{-1} y_0, \quad i = 1, \ldots\Big).$$
%\left(\prod_{i = 1}^k(A + s_i
%I_n)^{-1} \right)y_0$ $i = 1, \ldots, k$$. 
The shifts $s_i$ are selected by a
special procedure. There are different algorithms to compute the shifts 
(and the method proposed in this paper falls into this class, also there 
is a recent algorithm \cite{druskin-trks-2014} based on tangential interpolation). 
We use the RKSM method described in
\cite{druskin-rksm-2011} which has the publicly avaliable implementation.
The MATLAB code of both methods can be downloaded
 from \url{http://www.dm.unibo.it/~simoncin/software.html}.

Note that it is not fully fair to compare the efficiency 
of the ALR and KPIK
with RKSM. The first two methods use vectors from Krylov subspace and have
$2r + 1$ size of approximation subspaces at $r$ iteration step, 
in but ``pure'' RKSM method has a basis of size $(r+1)$ after $r$ iterations. 
We tried to ``extend'' the RKSM method by adding Krylov vectors. 
%It is not a trivial task due to a special functional, optimized in the 
%RKSM method to compute the shifts. However, we can mimic the 
%Krylov vectors by adding very large shift $s_i$at each odd iteration.
%This introduces an additional emprical parameter $s_i$ which have to be selected. 
%For our numerical examples we have found that the effect of such approach 
%is twofold. For small $r$, the convergence is improved. For larger $r$ 
%the method reaches a plateu. The exact reason for this behavior is to be investigated.
This extended approach we will call \emph{ERKSM}.  
\subsection{Efficient implementation}
\textcolor{black}
{
Note that the KPIK method has important feature. It uses only linear solver
with the matrix $A$ at every step of algorithm, so the factorization of the matrix $A$ can significantly
reduce the complexity. Another methods (ALR, RKSM, ERKSM) solve linear systems with different
matrices so this method is not applicable. 
}
\textcolor{black}
{
Instead, we propose to use the algebraic multigrid to speed-up the solution of linear systems 
 and to use the same multigrid hierarchy for different shifted systems. 
Algebraic multigrid (AMG) is a method for the solution of linear systems based on the multigrid principles, 
but requires no explicit knowledge of the problem geometry.
AMG determines coarse grids, intergrid
transfer operators, and coarse-grid equations
based solely on the matrix entries. 
Denote vector spaces $\mathbf{R^n}$ and $\mathbf{R^{n_c}}$ that
correspond to the fine and the coarse grids. Interpolation
(prolongation) maps the coarse grid to the
fine grid and is represented by the $n \times n_c$ matrix $P_c$ :
$\mathbf{R^n} \to \mathbf{R^{n_c}}$.
There are also few options for defining the coarse system $A_c$ and the most common
approach is to use the Galerkin projection $A_c = P_c^{\top} A P_c.$
}
\textcolor{black}
{
Our main idea is to use the algebraic multigrid method with the Galerkin
projection to construct the fast solver for ``shifted'' linear systems with matrices $(A
+ s I_n)$. Note that if the transfer operators are chosen to be the same for all shifts $s$, then
the coarse system matrix can be easy ``shifted'' $A_c + I_{n_c} = P_c^{\top}
(A +  s I_{n_c} ) P_c.$
}
\textcolor{black}
{
We use the Python implementation of AMG \cite{PyAMG-2}
to find the coarse grid hierarchy for the matrix $A$ once and then reuse this hierarchy for
shifted linear systems.  This trick is possible because
the matrix $A$ is stable and the shifts always have negative real parts. We found that
these approach works well both in symmetric and non-symmetric
cases in our numerical experiments.
}
%Although this approach may have such flaw that linear system solution timing increase with

\textcolor{black}
{
We have chosen $9$ methods for the comparison. First 4 methods are ALR,
KPIK, RKSM and ERKSM with shifted linear systems solved using direct solvers for sparse matrices. The next 4 methods 
use the modified AMG solver. And the last method is KPIK
with factorized inverse, and we denote it KPIK(LU). 
}
\subsection{Model problem 1: Laplace2D matrix}
\begin{figure}[H]
\scalebox{0.6}{\input{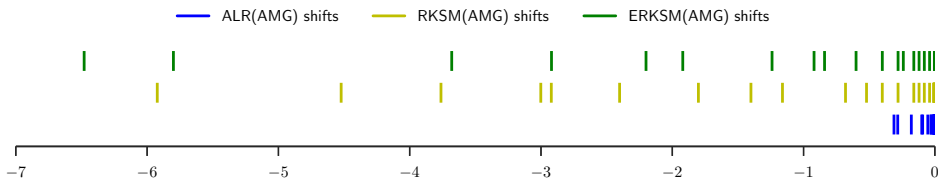}}
\caption{Shift distributions for different methods}\label{alr:shift-distrib}
\end{figure}

\textcolor{black}
{
We start our tests on symmetric problems. The first problem is the discretization of the two-dimensional
Laplace operator $$Lu = u_{xx} + u_{yy} $$ with Dirichlet boundary conditions on the unit square using 5-point stencil operator. 
The vector $y_0$ is obtained by the discretization on the grid of the
function $$f(x, y)= e^{-(x-0.5)^2 - 1.5(y-0.7)^2}.$$ The results are
shown in Table \ref{lyaptab:laplace2d}, final time shows the computational time in which the method reaches $10^{-8}$ 
accuracy.
}
\begin{table}[H]
\input{table_laplace2d.tex}
\end{table}
\textcolor{black}
{
Note that for the $256 \times 256$ grid the RKSM-based methods have significant precomputation time, in which 
the eigenvalue bounds are estimated. The KPIK(LU) method is significantly faster than other methods on all grids. 
Note that the AMG solver is slower than the direct solver as well. 
Another noticeable feature of the AMG solver is the significant variation of
average linear system solution time on the $256 \times 256$ grid for different shift selection strategies. 
The main reason of such variation is that the 
``shifted'' multigrid hierarchy may take more time in case of relatively big shifts. 
An example of the shift distribution is presented on Figure \ref{alr:shift-distrib}.
}
%\scalebox{0.13}{\input{laplace2d_shifts.pgf}}
%\input{laplace2d_shifts.pgf}
%\scalebox{0.5}{\includegraphics{laplace2d_shifts.pgf}}
%So, KPIK(AMG) method has the least average linear system time
%solution (0.41 secs), ALR(AMG) method has the second one (0.65 secs) and RKSM(AMG) and
%ERKSM(AMG)  methods have the biggest ones (0.91 and 1.04 correspondingly).
\subsection{Model problem 2: Laplace3D matrix}
\textcolor{black}
{
The second problem is taken from \cite[Example 5.3]{sim-krylov-2007}
and it is the discretization of the three-dimensional
Laplace operator  $$Lu = u_{xx} + u_{yy} + u_{zz}$$
with Dirichlet boundary conditions on the unit square using 7-point stencil operator. 
The vector $y_0$ was taken to be the vector of all ones. The results are
shown in Table \ref{lyaptab:sim53}, 
\begin{table}[h!]
    \input{table_sim53.tex}
\end{table}
%final time shows time at step in which method reaches $10^{-8}$ 
%accuracy.
}
%The KPIK(LU) method is faster only on the smallest grid. For larger grid size, 
%the KPIK(AMG) solver becomes faster. 
\subsection{Model problem 3}
\textcolor{black}
{
The third problem is taken from \cite[Example 5.1]{sim-krylov-2007} and describes a
model of heat flow with convection in the given domain. 
The associated differential operator is $$Lu = u_{xx} + u_{yy} - 10 x u_{x} - 1000 y u_{y}$$ on the unit square with Dirichlet boundary conditions.
%is the discretization of the two-dimensional
Matrices $A$ are obtained from the central finite difference
discretization of the differential operator using 5-point stencil
operator and are non-symmetric with complex eigenvalues.
Once again, the vector $y_0$ was
taken to be the vector of all ones.
The results are shown in Table \ref{lyaptab:sim51}.
}
%, final time shows time at step in which method reaches $10^{-8}$ 
%accuracy.
\begin{table}[h!]
\input{table_sim51.tex}\\ 
\end{table}
\textcolor{black}
{
For this problem KPIK method with LU factorization is the best
one on all grids.
}
\subsection{Model problem 4}
\textcolor{black}
{
The fourth problem is taken from \cite[Example 5.2]{sim-krylov-2007} and it is a three-dimensional variant of the previous problem.
The differential operator is given by
$$ Lu = u_{xx} + u_{yy} + u_{zz} - 10x u_{x} - 1000y u_{y} - u_{z}$$ with Dirichlet boundary conditions on the unit cube using 7-point stencil operator. 
Once again, the vector $y_0$ was taken to be the vector of all ones. 
The results are shown in Table \ref{lyaptab:sim52}.%
}
%, final time
%shows time at step in which method reaches $10^{-8}$ accuracy.
\begin{table}[h]
\input{table_sim52.tex}\\
\end{table}
\textcolor{black}
{
The LU-factorization takes too long on the finest grid for 3D problems so KPIK method which uses AMG solver is faster than
KPIK(LU) is faster. The ALR(AMG) method is the fastest for the $30 \times 30 \times 30$ grid.
%KPIK with
%LU-factorization is good at the coarsest grid and ALR with multigrid solver on the finest grid performs better than other methods.
}
\section{Conclusions and future work}
\textcolor{black}
{
The existing methods for the low-rank approximation to the solution of the Lyapunov equations perform 
quite well on different problems, provided the right implementation of the linear system solver with shifted matrices is available.
The ALR method in this paper produces different shifts in comparison with the RKSM solver: the shifts computed by the ALR method
are less spread and this gives a possibility to reuse the multigrid hierarchy in a more efficient way, and we think it is the most
promising feature of our method (and the number of iterations is similar to the RKSM/ERKSM methods).
}
\textcolor{black}
{
In our future work we plan to investigate the properties of the proposed approach. First of all, 
we would like to prove the convergence of the method. It will be also very interesting to generalize it  to the 
rank-$r$ right-hand side, and also to extend it to other matrix functions. The effect of inexact solvers on the convergence 
of the ALR method also need to be studied.
}
\section*{Acknowledgements}
 This work was supported by Russian Science Foundation grant 14-11-00659.
 We thank anonymous referees of the paper. Their comments helped to improve the quality of the paper, especially in the
 numerical experiments section.
We thank Dr. Vladimir Druskin for helpful comments on the RKSM approach
and the idea how to incorporate Krylov vectors into the existing RKSM code. 
We also thank Dr. Bart Vandreycken
for his comments on the initial draft of the paper.
\bibliographystyle{siam}
\bibliography{lyap}
\end{document}

%% file: table_laplace2d.tex
\begin{minipage}{\linewidth}
\centering
\captionof{table}{Model problem 1 timings}\label{lyaptab:laplace2d}\begin{tabular}{lllllr}
\toprule
grid &method &precomp. &1 solver time &final time &it-s/rank\\
\midrule
                       &ALR         & --         &0.0190      &0.2303      &10/21       \\
                       &ALR(AMG)    &0.0695      &0.0498      &0.5703      &10/21       \\
                       &KPIK(LU)    &\bf{0.0341} &\bf{0.0018} &\bf{0.1755} &15/31       \\
                       &KPIK(AMG)   &0.0695      &0.0319      &0.6104      &15/31       \\
64$\times$64           &KPIK        & --         &0.0196      &0.3654      &15/31       \\
                       &RKSM        &0.1345      &0.0194      &1.8460      &21/22       \\
                       &RKSM(AMG)   &0.2040      &0.0434      &2.4254      &21/22       \\
                       &ERKSM       &0.1345      &0.0210      &2.7598      &18/37       \\
                       &ERKSM(AMG)  &0.2040      &0.0533      &3.4273      &18/37       \\
\midrule
                       &ALR         & --         &0.0981      &1.3168      &12/25       \\
                       &ALR(AMG)    &0.1941      &0.1431      &2.5503      &12/25       \\
                       &KPIK(LU)    &\bf{0.1746} &\bf{0.0066} &\bf{1.0401} &20/41       \\
                       &KPIK(AMG)   &0.1941      &0.1018      &2.9015      &20/41       \\
128$\times$128         &KPIK        & --         &0.1051      &2.7162      &20/41       \\
                       &RKSM        &1.0598      &0.1024      &5.1593      &22/23       \\
                       &RKSM(AMG)   &1.2539      &0.1751      &7.4457      &23/24       \\
                       &ERKSM       &1.0598      &0.1063      &9.2354      &25/51       \\
                       &ERKSM(AMG)  &1.2539      &0.1819      &10.5334     &23/47       \\
\midrule
                       &ALR         & --         &0.6541      &10.7695     &15/31       \\
                       &ALR(AMG)    &\bf{0.6764} &0.6083      &13.1072     &15/31       \\
                       &KPIK(LU)    &1.1841      &\bf{0.0399} &\bf{9.5395} &26/53       \\
                       &KPIK(AMG)   &\bf{0.6764} &0.4087      &17.2390     &26/53       \\
256$\times$256         &KPIK        & --         &0.6728      &21.4517     &26/53       \\
                       &RKSM        &29.7068     &0.7498      &55.8941     &27/28       \\
                       &RKSM(AMG)   &30.3832     &0.9082      &59.5665     &26/27       \\
                       &ERKSM       &29.7068     &0.7158      &58.3618     &24/49       \\
                       &ERKSM(AMG)  &30.3832     &1.0387      &68.0542     &24/49       \\
\bottomrule
\end{tabular}\par
\bigskip
\end{minipage}

%% file: table_sim53.tex
\begin{minipage}{\linewidth}
\centering
\captionof{table}{Model problem 2 timings}\label{lyaptab:sim53}\begin{tabular}{lllllr}
\toprule
grid &method &precomp. &1 solver time &final time &it-s/rank\\
\midrule
                       &ALR         & --         &0.0147      &0.0822      &5/11        \\
                       &ALR(AMG)    &0.0712      &0.0113      &0.1286      &5/11        \\
                       &KPIK(LU)    &0.0163      &\bf{0.0006} &\bf{0.0303} &6/13        \\
                       &KPIK(AMG)   &0.0712      &0.0136      &0.1521      &6/13        \\
10$\times$10$\times$10 &KPIK        & --         &0.0116      &0.0684      &6/13        \\
                       &RKSM        &\bf{0.0054} &0.0102      &0.3186      &9/10        \\
                       &RKSM(AMG)   &0.0766      &0.0099      &0.3205      &9/10        \\
                       &ERKSM       &\bf{0.0054} &0.0142      &0.3081      &6/13        \\
                       &ERKSM(AMG)  &0.0766      &0.0146      &0.3729      &6/13        \\
\midrule
                       &ALR         & --         &0.6237      &3.8685      &7/15        \\
                       &ALR(AMG)    &0.1385      &0.0993      &0.7724      &7/15        \\
                       &KPIK(LU)    &0.8283      &\bf{0.0085} &0.9671      &8/17        \\
                       &KPIK(AMG)   &0.1385      &0.0766      &\bf{0.7419} &8/17        \\
20$\times$20$\times$20 &KPIK        & --         &0.6313      &4.4834      &8/17        \\
                       &RKSM        &\bf{0.0608} &0.6303      &6.0856      &10/11       \\
                       &RKSM(AMG)   &0.1993      &0.0846      &1.3038      &10/11       \\
                       &ERKSM       &\bf{0.0608} &0.6811      &6.9582      &10/21       \\
                       &ERKSM(AMG)  &0.1993      &0.0832      &1.5275      &9/19        \\
\midrule
                       &ALR         & --         &8.9842      &62.6705     &8/17        \\
                       &ALR(AMG)    &0.4265      &0.3449      &3.0680      &8/17        \\
                       &KPIK(LU)    &19.0281     &\bf{0.0535} &19.9242     &10/21       \\
                       &KPIK(AMG)   &0.4265      &0.2354      &\bf{2.9497} &10/21       \\
30$\times$30$\times$30 &KPIK        & --         &9.1710      &82.3441     &10/21       \\
                       &RKSM        &\bf{0.3609} &9.0977      &120.1328    &14/15       \\
                       &RKSM(AMG)   &0.7874      &0.3443      &6.2245      &14/15       \\
                       &ERKSM       &\bf{0.3609} &10.2733     &115.8494    &12/25       \\
                       &ERKSM(AMG)  &0.7874      &0.3836      &6.5964      &11/23       \\
\bottomrule
\end{tabular}\par
\bigskip
\end{minipage}

%% file: table_sim51.tex
\begin{minipage}{\linewidth}
\centering
\captionof{table}{Model problem 3 timings}\label{lyaptab:sim51}\begin{tabular}{lllllr}
\toprule
grid &method &precomp. &1 solver time &final time &it-s/rank\\
\midrule
                       &ALR         & --         &0.0184      &0.1834      &8/17        \\
                       &ALR(AMG)    &0.1744      &0.0253      &0.3649      &8/17        \\
                       &KPIK(LU)    &\bf{0.0197} &\bf{0.0017} &\bf{0.1575} &11/23       \\
                       &KPIK(AMG)   &0.1744      &0.0284      &0.5586      &11/23       \\
64$\times$64           &KPIK        & --         &0.0250      &0.3809      &11/23       \\
                       &RKSM        &0.0452      &0.0177      &0.7313      &12/13       \\
                       &RKSM(AMG)   &0.2196      &0.0185      &0.7349      &12/13       \\
                       &ERKSM       &0.0452      &0.0183      &0.7154      &9/19        \\
                       &ERKSM(AMG)  &0.2196      &0.0187      &1.0360      &10/21       \\
\midrule
                       &ALR         & --         &0.1189      &1.2058      &10/21       \\
                       &ALR(AMG)    &0.5329      &0.0875      &1.6560      &10/21       \\
                       &KPIK(LU)    &\bf{0.1136} &\bf{0.0098} &\bf{0.6804} &11/23       \\
                       &KPIK(AMG)   &0.5329      &0.1052      &2.0300      &11/23       \\
128$\times$128         &KPIK        & --         &0.1278      &1.5739      &11/23       \\
                       &RKSM        &0.6967      &0.1241      &3.0168      &13/14       \\
                       &RKSM(AMG)   &1.2296      &0.0758      &2.7024      &13/14       \\
                       &ERKSM       &0.6967      &0.1225      &3.5379      &12/25       \\
                       &ERKSM(AMG)  &1.2296      &0.0698      &3.2370      &11/23       \\
\midrule
                       &ALR         & --         &1.1270      &14.2010     &10/21       \\
                       &ALR(AMG)    &\bf{0.7507} &0.3281      &5.1982      &10/21       \\
                       &KPIK(LU)    &1.0185      &\bf{0.0313} &\bf{3.0908} &12/25       \\
                       &KPIK(AMG)   &\bf{0.7507} &0.2982      &5.6122      &12/25       \\
256$\times$256         &KPIK        & --         &1.4760      &18.2555     &12/25       \\
                       &RKSM        &6.8556      &1.0914      &27.0884     &17/18       \\
                       &RKSM(AMG)   &7.6063      &0.2696      &14.1912     &17/18       \\
                       &ERKSM       &6.8556      &1.1761      &29.3055     &15/31       \\
                       &ERKSM(AMG)  &7.6063      &0.2524      &15.1329     &13/27       \\
\bottomrule
\end{tabular}\par
\bigskip
\end{minipage}

%% file: table_sim52.tex
\begin{minipage}{\linewidth}
\centering
\captionof{table}{Model problem 4 timings}\label{lyaptab:sim52}\begin{tabular}{lllllr}
\toprule
grid &method &precomp. &1 solver time &final time &it-s/rank\\
\midrule
                       &ALR         & --         &0.0117      &0.0671      &5/11        \\
                       &ALR(AMG)    &0.0444      &0.0066      &0.0952      &5/11        \\
                       &KPIK(LU)    &0.0181      &\bf{0.0004} &\bf{0.0386} &7/15        \\
                       &KPIK(AMG)   &0.0444      &0.0099      &0.1206      &7/15        \\
10$\times$10$\times$10 &KPIK        & --         &0.0257      &0.1760      &7/15        \\
                       &RKSM        &\bf{0.0041} &0.0108      &0.2189      &7/8         \\
                       &RKSM(AMG)   &0.0485      &0.0094      &0.2259      &7/8         \\
                       &ERKSM       &\bf{0.0041} &0.0128      &0.2201      &5/11        \\
                       &ERKSM(AMG)  &0.0485      &0.0085      &0.2276      &5/11        \\
\midrule
                       &ALR         & --         &0.6424      &3.2881      &6/13        \\
                       &ALR(AMG)    &0.2231      &0.0466      &\bf{0.5020} &6/13        \\
                       &KPIK(LU)    &0.6687      &\bf{0.0086} &0.8351      &9/19        \\
                       &KPIK(AMG)   &0.2231      &0.0436      &0.6787      &9/19        \\
20$\times$20$\times$20 &KPIK        & --         &0.8533      &7.0838      &9/19        \\
                       &RKSM        &\bf{0.0265} &0.6610      &5.5919      &9/10        \\
                       &RKSM(AMG)   &0.2496      &0.0491      &0.9061      &9/10        \\
                       &ERKSM       &\bf{0.0265} &0.7303      &4.8876      &7/15        \\
                       &ERKSM(AMG)  &0.2496      &0.0541      &0.9108      &6/13        \\
\midrule
                       &ALR         & --         &9.2135      &53.2061     &7/15        \\
                       &ALR(AMG)    &0.9827      &0.1987      &\bf{2.3924} &7/15        \\
                       &KPIK(LU)    &10.1743     &\bf{0.0497} &11.0392     &9/19        \\
                       &KPIK(AMG)   &0.9827      &0.2055      &3.0282      &9/19        \\
30$\times$30$\times$30 &KPIK        & --         &12.7750     &105.1641    &9/19        \\
                       &RKSM        &\bf{0.1215} &9.3640      &85.8903     &10/11       \\
                       &RKSM(AMG)   &1.1042      &0.1792      &3.2822      &10/11       \\
                       &ERKSM       &\bf{0.1215} &9.5372      &58.2637     &7/15        \\
                       &ERKSM(AMG)  &1.1042      &0.1609      &3.2703      &8/17        \\
\bottomrule
\end{tabular}\par
\bigskip
\end{minipage}